\documentclass{amsart}
\usepackage{times,amssymb,amscd,graphics,axodraw,verbatim}
\usepackage[stdtext,vcentermath]{myyoungtab}
\usepackage[all]{xy}

\numberwithin{equation}{section}

\newtheorem{prop}{Proposition}
\newtheorem{theorem}[prop]{Theorem}
\newtheorem{corollary}[prop]{Corollary}
\newtheorem{lemma}[prop]{Lemma}
\newtheorem{conjecture}[prop]{Conjecture}

\theoremstyle{definition}
\newtheorem{definition}[prop]{Definition}
\newtheorem{example}[prop]{Example}
\newtheorem{remark}[prop]{Remark}

\numberwithin{prop}{section}

\newcommand{\A}{\mathcal{A}}
\newcommand{\cc}{cc}
\newcommand{\Conf}{\mathrm{C}}
\newcommand{\Confb}{\overline{\Conf}}

\newcommand{\geh}{\mathfrak{g}}
\newcommand{\HH}{\mathcal{H}}
\newcommand{\inner}[2]{\langle #1\,,\,#2\rangle}
\newcommand{\J}{\overline{I}}
\newcommand{\Jb}{\overline{J}}
\newcommand{\la}{\lambda}
\newcommand{\La}{\Lambda}

\newcommand{\Mb}{\overline{M}}
\newcommand{\N}{\mathbb{N}}
\newcommand{\nub}{\overline{\nu}}
\newcommand{\Path}{\mathcal{P}}
\newcommand{\Pathb}{\overline{\Path}}
\newcommand{\Phib}{\overline{\Phi}}
\newcommand{\pr}{\mathrm{pr}}
\newcommand{\qbin}[2]{\genfrac{[}{]}{0pt}{}{#1}{#2}}
\newcommand{\RC}{\mathrm{RC}}
\newcommand{\RCb}{\overline{\RC}}
\newcommand{\RCt}{\widetilde{\RC}}

\newcommand{\SA}{\mathcal{SA}}
\newcommand{\ve}{\varepsilon}
\newcommand{\vp}{\varphi}
\newcommand{\wt}{\mathrm{wt}}
\newcommand{\Xb}{\overline{X}}
\newcommand{\Z}{\mathbb{Z}}

\begin{document}

\title{Crystal structure on rigged configurations}

\author[A.~Schilling]{Anne Schilling}
\address{Department of Mathematics, University of California, One Shields
Avenue, Davis, CA 95616-8633, U.S.A.}
\email{anne@math.ucdavis.edu}
\urladdr{http://www.math.ucdavis.edu/\~{}anne}
\thanks{\textit{Date:} July 2005}
\thanks{Partially supported by NSF grants DMS-0200774 and DMS-0501101.}
 
\begin{abstract}
Rigged configurations are combinatorial objects originating from the Bethe
Ansatz, that label highest weight crystal elements. In this paper a new
\textit{unrestricted} set of rigged configurations is introduced for types $ADE$ 
by constructing a crystal structure on the set of rigged configurations. 
In type $A$ an explicit characterization of unrestricted rigged configurations 
is provided which leads to a new fermionic formula for unrestricted Kostka
polynomials or $q$-supernomial coefficients. The affine crystal structure for 
type $A$ is obtained as well.
\end{abstract}

\maketitle

\section{Introduction}
There are (at least) two main approaches to solvable lattice models
and their associated quantum spin chains: the Bethe Ansatz~\cite{Bethe:1931} 
and the corner transfer matrix method~\cite{Baxter:1982}.

In his 1931 paper~\cite{Bethe:1931}, Bethe solved the Heisenberg spin chain
based on the string hypothesis which asserts that the eigenvalues of
the Hamiltonian form certain strings in the complex plane as the size
of the system tends to infinity. The Bethe Ansatz has been applied to many
models to prove completeness of the Bethe vectors.
The eigenvalues and eigenvectors of the Hamiltonian are 
indexed by rigged configurations. However, numerical studies indicate 
that the string hypothesis is not always true~\cite{ADM:1992}.

The corner transfer matrix (CTM) method, introduced by Baxter~\cite{Baxter:1982},
labels the eigenvectors by one-dimensional lattice paths. These lattice paths have 
a natural interpretation in terms of Kashiwara's crystal base 
theory~\cite{Kash:1990,Kash:1994}, namely as highest weight crystal elements in a 
tensor product of finite-dimensional crystals.

\begin{figure}
\resizebox{10cm}{!}{\includegraphics{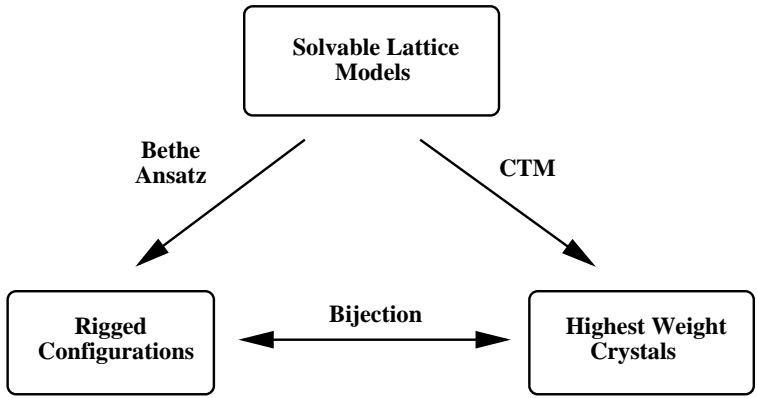}}
\caption{Schematic origin of rigged configurations and crystal paths
\label{fig:scheme}}
\end{figure}

Even though neither the Bethe Ansatz nor the corner transfer matrix
method are mathematically rigorous, they suggest that the existence of
a bijection between the two index sets, namely rigged configurations on the
one hand and highest weight crystal paths on the other (see Figure~\ref{fig:scheme}).
For the special case when the spin chain is defined on 
$V_{(\mu_1)}\otimes V_{(\mu_2)} \otimes \cdots \otimes V_{(\mu_k)}$, where $V_{(\mu_i)}$
is the irreducible $\mathrm{GL}(n)$ representation indexed by the partition $(\mu_i)$ 
for $\mu_i\in \N$, a bijection between rigged configurations and semi-standard
Young tableaux was given by Kerov, Kirillov and Reshetikhin~\cite{KKR:1986,KR:1988}.
This bijection was proven and extended to the case when the $(\mu_i)$ are any sequence of
rectangles in~\cite{KSS:2002}. The bijection has many amazing properties. For
example it takes the cocharge statistics $\cc$ defined on rigged configurations
to the coenergy statistics $D$ defined on crystals.

Rigged configurations and crystal paths also exist for other types.
In~\cite{HKOTT:2001,HKOTY:1999} the existence of Kirillov--Reshetikhin
crystals $B^{r,s}$ was conjectured, which can be naturally associated with
the dominant weight $s\La_r$ where $s$ is a positive integer and $\La_r$ is the
$r$-th fundamental weight of the underlying algebra of finite type. For a tensor product 
of Kirillov--Reshetikhin crystals $B=B^{r_k,s_k}\otimes \cdots \otimes B^{r_1,s_1}$ 
and a dominant weight $\La$ let $\Pathb(B,\La)$ be the set of all highest weight 
elements of weight $\La$ in $B$. In the same papers~\cite{HKOTT:2001,HKOTY:1999}, 
fermionic formulas $\Mb(L,\La)$ for the one-dimensional configuration sums 
$\Xb(B,\La):=\sum_{b\in\Pathb(B,\La)} q^{D(b)}$ were conjectured. The fermionic formulas 
admit a combinatorial interpretation in terms of the set of rigged configurations 
$\RCb(L,\La)$, where $L$ is the multiplicity array of $B$ (see Section~\ref{sec:def rc}). 
A statistic preserving bijection $\Phi:\Pathb(B,\La)\to\RCb(L,\La)$ has been proven in 
various cases~\cite{KSS:2002,OSS:2003,S:2005,SS:2005} which implies the following identity
\begin{equation}\label{eq:XM intro}
\Xb(B,\La):=\sum_{b\in\Pathb(B,\La)} q^{D(b)} 
= \sum_{(\nu,J)\in\RCb(L,\La)} q^{\cc(\nu,J)}=:\Mb(L,\La).
\end{equation}
Since the sets in~\eqref{eq:XM intro} are finite, these are polynomials in $q$.
When $B=B^{1,s_k}\otimes \cdots \otimes B^{1,s_1}$ of type $A$, they are 
none other than the Kostka--Foulkes polynomials.

Rigged configurations corresponding to highest weight crystal paths are only the 
tip of an iceberg. In this paper we extend the definition of rigged configurations 
to all crystal elements in types $ADE$ by the explicit construction of a crystal 
structure on the set of \textit{unrestricted} rigged configurations 
(see Definition~\ref{def:crystal}). The proof uses Stembridge's local characterization 
of simply-laced crystals~\cite{St:2003}. For nonsimply-laced algebras, the local rules 
provided in~\cite{St:2003} are still necessary, but no longer sufficient 
conditions to characterize crystals. Crystal operators for rigged configurations 
associated to nonsimply-laced algebras can be constructed from the ones presented here
via ``folding'' of the Dynkin diagrams as in the construction of virtual 
crystals~\cite{OSS:2003a,OSS:2003b}.

The equivalence of the crystal structures on rigged configurations and crystal paths 
together with the correspondence for highest weight vectors yields the equality 
of generating functions in analogy to~\eqref{eq:XM intro} (see Theorem~\ref{thm:bij new}
and Corollary~\ref{cor:X=M}). Denote the unrestricted
set of paths and rigged configurations by $\Path(B,\La)$ and $\RC(L,\La)$, respectively.
The corresponding generating functions $X(B,\La)=M(L,\La)$ are unrestricted generalized 
Kostka polynomials or $q$-supernomial coefficients.
A direct bijection $\Phi:\Path(B,\La)\to \RC(L,\La)$ for type $A$ along the lines 
of~\cite{KSS:2002} is constructed in~\cite{DS:2004,DS:2005}.

Rigged configurations are closely tied to fermionic formulas.
Fermionic formulas are explicit expressions for the partition function
of the underlying physical model which reflect their particle structure.
For more details regarding the background of fermionic formulas
see~\cite{HKOTT:2001,KKMM:1993a,KKMM:1993b}. For type $A$ we obtain
an explicit characterization of the unrestricted rigged configurations
in terms of lower bounds on quantum numbers (see Definition~\ref{def:extended}
and Theorem~\ref{thm:ext=unres}) which yields a new fermionic formula
for unrestricted Kostka polynomials of type $A$ (see Equation~\eqref{eq:fermi}).
Surprisingly, this formula is different from the fermionic formulas 
in~\cite{HKKOTY:1999,Kir:2000} obtained in the special cases of 
$B=B^{1,s_k}\otimes \cdots \otimes B^{1,s_1}$ and $B=B^{r_k,1}\otimes \cdots 
\otimes B^{r_1,1}$. The rigged configurations corresponding to the fermionic formulas 
of~\cite{HKKOTY:1999,Kir:2000} were related to ribbon tableaux and the cospin generating 
functions of Lascoux, Leclerc, Thibon~\cite{LLT:1997,LT:2000} in reference~\cite{S:2002}.
To distinguish these rigged configurations from the ones introduced in this
paper, let us call them ribbon rigged configurations.

The Lascoux--Leclerc--Thibon (LLT) polynomials~\cite{LLT:1997,LT:2000}
have recently made their debut in the theory of Macdonald polynomials
in the seminal paper by Haiman, Haglund, Loehr~\cite{HHL:2005}. The main
obstacle in obtaining a combinatorial formula for the Macdonald--Kostka
polynomials is the Schur positivity of certain LLT polynomials.
A related problem is the conjecture of Kirillov and Shimozono~\cite{KS:2002}
that the cospin generating function of ribbon tableaux equals the generalized
Kostka polynomial. A possible avenue to prove this conjecture would be a
direct bijection between the unrestricted rigged configurations of this
paper and ribbon rigged configurations.

For type $A$ we can also describe the affine crystal operators $e_0$ and
$f_0$ on rigged configurations (see Section~\ref{sec:affine}). A level-$\ell$
restricted element $b$ in a crystal $B$ is characterized by $e_0^{\ell+1}b=0$.
It is striking that the description of the unrestricted rigged 
configurations of type $A$ (see Definition~\ref{def:extended})
is very similar to the characterization of level-restricted rigged 
configurations as given in~\cite{SS:2001}. Instead of a modification of
the vacancy numbers, which comprise upper bounds for the quantum numbers,
an alteration of the lower bounds occurs. In both cases the modification
is governed by a set of tableaux depending on the weight $\La$. It would be interesting
to understand the relation between the conditions imposed by classical restriction
and level-restriction in a more precise manner.

One of the motivations for considering unrestricted rigged configurations
was Takagi's work~\cite{T:2005} on the inverse scattering transform, which
provides a bijection between states in the $\mathfrak{sl}_2$ box ball system
and rigged configurations. In this setting rigged configurations play the role
of action-angle variables. Box ball systems can be produced from crystals
of solvable lattice models for algebras other than 
$\mathfrak{sl}_2$~\cite{HIK:1999,HKT:2000,HHIKTT:2001}.
The inverse scattering transform can be generalized to the $\mathfrak{sl}_n$
case~\cite{KOSTY:2006}, which should give a box-ball interpretation
of the unrestricted rigged configurations presented here.

Another motivation for the study of unrestricted configuration sums,
fermionic formulas and associated rigged configurations is their appearance
in generalizations of the Bailey lemma~\cite{ASW:1999,W:2004}.
The Andrews--Bailey construction~\cite{Andrews:1984,Bailey:1949} relies
on an iterative transformation property of the $q$-binomial coefficient,
which is one of the simplest unrestricted configuration sums, and can be
used to prove infinite families of Rogers--Ramanujan type identities.
The explicit formulas provided in this paper might trigger further progress
towards generalizations to higher-rank or other types of the 
Andrews--Bailey construction.

The paper is organized as follows. In Section~\ref{sec:axiom} we review
basics about crystal bases and Stembridge's local characterization of
crystals. In Section~\ref{sec:rc} we define rigged configurations
and the new crystal structure for types $ADE$. Section~\ref{sec:type A}
is devoted to type $A$, where we give an explicit characterization
of the unrestricted rigged configurations in Section~\ref{sec:lower bound},
a new fermionic formula for unrestricted Kostka polynomials in
Section~\ref{sec:fermi}, and the affine crystal structure in 
Section~\ref{sec:affine}.

\subsection*{Acknowledgment} I would like to thank Mark Haiman, Mark Shimozono 
and John Stembridge for helpful discussions, Peter Littelmann for drawing my 
attention to reference~\cite{St:2003}, and Masato Okado for his explanations 
of the box ball system~\cite{KOSTY:2006}. I would also like to thank Lipika 
Deka for collaboration on~\cite{DS:2004,DS:2005}.

\section{Crystal graphs}

We review the axiomatic definition of crystal graphs in Section~\ref{sec:axiom}
and the local characterization of crystals corresponding to representations of 
simply-laced algebras provided by Stembridge~\cite{St:2003} in Section~\ref{sec:char}.
In Section~\ref{sec:KR} we review the main properties of Kirillov--Reshetikhin
crystals.

\subsection{Axiomatic definition} \label{sec:axiom}

Kashiwara~\cite{Kash:1990,Kash:1994} introduced a \textbf{crystal}
as an edge-colored directed graph satisfying a simple set of axioms.
Let $\geh$ be a symmetrizable Kac--Moody algebra with associated root, coroot
and weight lattices $Q,Q^\vee,P$. Let $I$ be the index set of the Dynkin 
diagram and denote the simple roots, simple coroots and fundamental weights 
by $\alpha_i$, $h_i$ and $\Lambda_i$ ($i\in I$), respectively.
There is a natural pairing $\inner{\cdot}{\cdot}:Q^\vee\otimes P \rightarrow
\Z$ defined by $\inner{h_i}{\La_j}=\delta_{ij}$.

The vertices of the crystal graph are elements of a set $B$.
The edges of the crystal graph are colored by the index set
$I$. A $P$-weighted $I$-crystal satisfies the following properties:
\begin{enumerate}
\item Fix an $i\in I$. If all edges are removed except those colored $i$, the
connected components are finite directed linear paths called the
\textbf{$i$-strings} of $B$. Given $b\in B$, define $f_i(b)$
(resp. $e_i(b)$) to be the vertex following (resp. preceding) $b$
in its $i$-string; if there is no such vertex, declare $f_i(b)$ (resp. $e_i(b)$)
to be undefined. Define $\vp_i(b)$ (resp. $\ve_i(b)$) to be the number of arrows 
from $b$ to the end (resp. beginning) of its $i$-string.
\item There is a function $\wt:B\rightarrow P$ such that
\begin{equation*}
\begin{split}
\wt(f_i(b))&=\wt(b)-\alpha_i \\
\vp_i(b)-\ve_i(b) &= \inner{h_i}{\wt(b)}.
\end{split}
\end{equation*}
\end{enumerate}

\subsection{Local characterization of crystals} \label{sec:char}

Kashiwara~\cite{Kash:1990,Kash:1994} constructed crystal
graphs for representations of $U_q(\geh)$. Crystals of representations
form a special subclass of the set of all crystals.
In~\cite{St:2003}, Stembridge determines a simple set of local
axioms that uniquely characterize the crystals corresponding to
representations for simply-laced algebras.

Let $A=[A_{ij}]_{i,j\in I}$ be the Cartan matrix of a simply-laced
Kac--Moody algebra $\geh$. Let $X$ be an edge-colored graph.
Stembridge~\cite{St:2003} introduces the notion of $A$-regularity
by requiring the conditions (P1)-(P6), (P5'), (P6') to hold.
\begin{enumerate}
\item[(P1)] All monochromatic directed paths in $X$ have finite length.
In particular $X$ has no monochromatic circuits.
\item[(P2)] For every vertex $x$ and every $i\in I$, there is at most
one edge $y \stackrel{i}{\longrightarrow} x$, and dually, at most one edge
$x\stackrel{i}{\longrightarrow} z$.
\end{enumerate}

In the notation of the previous section, the relation $f_i(x)=y$ or
equivalently $e_i(y)=x$ is graphically depicted by 
$x\stackrel{i}{\longrightarrow} y$. Set $\delta_i(x)=-\ve_i(x)$ with
$\ve_i(x)$ and $\vp_i(x)$ as defined in Section~\ref{sec:axiom}. Define
\begin{equation*}
\Delta_i \delta_j(x)=\delta_j(e_i x)-\delta_j(x), \qquad
\Delta_i \vp_j(x)=\vp_j(e_i x)-\vp_j(x),
\end{equation*}
whenever $e_i x$ is defined, and
\begin{equation*}
\nabla_i \delta_j(x)=\delta_j(x)-\delta_j(f_i x), \qquad
\nabla_i \vp_j(x)=\vp_j(x)-\vp_j(f_i x),
\end{equation*}
whenever $f_i x$ is defined.

For fixed $x\in X$ and a distinct pair $i,j\in I$, assuming that $e_i x$ is 
defined, require
\begin{enumerate}
\item[(P3)] $\Delta_i \delta_j(x) + \Delta_i \vp_j(x) = A_{ij}$, and
\item[(P4)] $\Delta_i \delta_j(x)\le 0$, $\Delta_i \vp_j(x)\le 0$.
\end{enumerate}
Note that for simply-laced algebras $A_{ij}\in \{0,-1\}$ for $i,j\in I$ distinct.
Hence (P3) and (P4) allow for only three possibilities:
\begin{equation*}
(A_{ij},\Delta_i\delta_j(x),\Delta_i\vp_j(x))=(0,0,0),(-1,-1,0),(-1,0,-1).
\end{equation*} 

Assuming that $e_ix$ and $e_jx$ both exist, we require
\begin{enumerate}
\item[(P5)] $\Delta_i\delta_j(x)=0$ implies $y:=e_ie_jx=e_je_ix$ and
$\nabla_j\vp_i(y)=0$.
\item[(P6)] $\Delta_i\delta_j(x)=\Delta_j\delta_i(x)=-1$ implies
$y:=e_i e_j^2 e_i x= e_j e_i^2 e_j x$ and $\nabla_i\vp_j(y)=\nabla_j\vp_i(y)=-1$.
\end{enumerate}

Dually, assuming that $f_ix$ and $f_jx$ both exist, we require
\begin{enumerate}
\item[(P5')] $\nabla_i\vp_j(x)=0$ implies $y:=f_i f_j x = f_j f_i x$ and
$\Delta_j\delta_i(y)=0$.
\item[(P6')] $\nabla_i\vp_j=\nabla_j\vp_i(x)=-1$ implies
$y:=f_i f_j^2 f_i x= f_j f_i^2 f_j x$ and $\Delta_i\delta_j(y)=\Delta_j
\delta_i(y)=-1$.
\end{enumerate}

\begin{definition} \cite[Definition 1.1]{St:2003}
Let $A$ be a simply-laced Cartan matrix. An edge-colored directed graph is
\textbf{$A$-regular} if it satisfies (P1)-(P6) and (P5')-(P6').
\end{definition}

Stembridge proved~\cite[Proposition 1.4]{St:2003} that any two $A$-regular
posets $P,P'$ with maximal elements $x,x'$ are isomorphic if and only if
$\vp_i(x)=\vp_i(x')$ for all $i\in I$. Moreover this isomorphism is unique.
Let $\La=\sum_{i\in I} \mu_i \La_i$ be a dominant weight. 
Denote by $B(\La)$ the unique $A$-regular poset with maximal element $b$ such that 
$\vp_i(b)=\mu_i$ for all $i\in I$.

\begin{theorem} \cite[Theorem 3.3]{St:2003} \label{thm:stembridge}
If $A$ is a simply-laced Cartan matrix, then the crystal graph of the 
irreducible $U_q(A)$-module of highest weight $\La$ is $B(\La)$.
\end{theorem}

\subsection{Kirillov--Reshetikhin crystals} \label{sec:KR}

Kirillov--Reshetikhin crystals are crystals for finite-dimensional irreducible
modules over quantum affine algebras. The irreducible finite-dimensional
$U'_q(\geh)$-modules were classified by Chari and Pressley~\cite{CP:1995,CP:1998}
in terms of Drinfeld polynomials. Here $U'_q(\geh)$ is the derived algebra without
the generator $q^d$, where $d$ is degree operator in $\geh$.
The Kirillov--Reshetikhin modules $W^{r,s}$, labeled by a Dynkin node $r$ of the 
underlying algebra of finite type and a positive integer $s$, form a special class of these
finite-dimensional modules.  They naturally correspond to the weight
$s\Lambda_r$, where $\Lambda_r$ is the $r$-th fundamental weight of
$\geh$. It was conjectured in~\cite{HKOTT:2001,HKOTY:1999}, that there exists a crystal
$B^{r,s}$ for each $W^{r,s}$. In general, the existence of $B^{r,s}$ is still an
open question. For type $A_n^{(1)}$ the crystal $B^{r,s}$ is known to
exist~\cite{KKMMNN:1992} and its combinatorial structure has been studied~\cite{Sh:2002}.
For other types, the existence and combinatorial structure of $B^{r,s}$ has been
considered in various articles (see for example~\cite{K:1999,KKM:1994,KKMMNN:1992}).
As classical crystals the Kirillov--Reshetikhin crystals are isomorphic to
\begin{equation*}
B^{r,s} \simeq B(s\La_r) \oplus \bigoplus_{\La} B(\La),
\end{equation*}
where $B(\La)$ is the classically highest weight crystal of highest weight $\La$
and the sum is over a particular set of weights contained in $s\La_r$
(for more details see~\cite{HKOTT:2001}).

\section{Crystal structure on rigged configurations} \label{sec:rc}

In this section we define a crystal structure on rigged configurations.
As alluded to in the introduction rigged configurations form a combinatorial set to 
index the eigenvalues and eigenvectors of the Hamiltonian of an exactly solvable lattice 
model. The simplest version of rigged configurations appeared in Bethe's original
paper~\cite{Bethe:1931} and was later generalized by Kerov, Kirillov and 
Reshetikhin~\cite{KKR:1986,KR:1988} to models with $\mathrm{GL}(n)$ symmetry.
Since the eigenvectors of the Hamiltonian can also be viewed as highest
weight vectors, one expects a bijection between rigged configurations
and semi-standard Young tableaux in the $\mathrm{GL}(n)$ case. Such a bijection
was given in~\cite{KR:1988,KSS:2002}. Rigged configurations for other types
follow from the fermionic formulas given in~\cite{HKOTT:2001,HKOTY:1999}
and they correspond to highest weight crystal elements~\cite{OSS:2003,S:2005,SS:2005}.
Here we extend the notion of rigged configurations to
non-highest weight elements called unrestricted rigged configurations and define 
a crystal structure on this set. In Section~\ref{sec:def rc} we review the
definition and known results about the usual rigged configurations.
In Section~\ref{sec:def rc crystal} the set of unrestricted rigged configurations
is introduced and the crystal structure is defined for types $ADE$ 
(Definition~\ref{def:crystal}).
This leads to a bijection between crystal paths and unrestricted rigged
configurations (Theorem~\ref{thm:bij new}) and the equality of generating
functions (Corollary~\ref{cor:X=M}).

\subsection{Definition of rigged configurations}\label{sec:def rc}
Let $\geh$ be a simple simply-laced affine Kac--Moody algebra.
Define $\J=I\setminus\{0\}$ the index set of the underlying algebra of finite
type and set $\HH=\J \times \Z_{>0}$. The (highest-weight) rigged configurations
are indexed by a multiplicity array $L=(L_i^{(a)}\mid (a,i)\in \HH)$ of
nonnegative integers and a dominant weight $\La$.
The sequence of partitions $\nu=\{\nu^{(a)}\mid a\in \J \}$ is a
\textbf{$(L,\La)$-configuration} if
\begin{equation}\label{eq:conf}
\sum_{(a,i)\in\HH} i m_i^{(a)} \alpha_a = \sum_{(a,i)\in\HH} i
L_i^{(a)} \La_a- \La,
\end{equation}
where $m_i^{(a)}$ is the number of parts of length $i$ in partition
$\nu^{(a)}$. Denote the set of all $(L,\La)$-configurations by $\Conf(L,\La)$.
The \textbf{vacancy number} of a configuration is defined as
\begin{equation}\label{eq:vac}
p_i^{(a)}=\sum_{j\ge 1} \min(i,j) L_j^{(a)}
 - \sum_{(b,j)\in \HH} (\alpha_a | \alpha_b) \min(i,j)m_j^{(b)}.
\end{equation}
Here $(\cdot | \cdot )$ is the normalized invariant form on the weight lattice $P$
such that $A_{ab}=(\alpha_a | \alpha_b)$ is the Cartan matrix.
The $(L,\La)$-configuration $\nu$ is \textbf{admissible} if $p^{(a)}_i\ge 0$ for 
all $(a,i)\in\HH$, and the set of admissible $(L,\La)$-configurations is denoted
by $\Confb(L,\La)$.

A rigged configuration is an admissible configuration together with
a set of labels of quantum numbers. A partition can be viewed as
a multiset of positive integers.
A rigged partition is by definition a finite multiset of
pairs $(i,x)$ where $i$ is a positive integer and
$x$ is a nonnegative integer.  The pairs $(i,x)$ are referred to
as strings; $i$ is referred to as the
length or size of the string and $x$ as the \textbf{label} or
\textbf{quantum number} of the string.  A rigged partition is
said to be a rigging of the partition $\rho$ if
the multiset, consisting of the sizes of the strings,
is the partition $\rho$.  So a rigging of $\rho$
is a labeling of the parts of $\rho$ by nonnegative integers,
where one identifies labelings that differ only by
permuting labels among equal sized parts of $\rho$.

A rigging $J$ of the $(L,\La)$-configuration $\nu$ is a sequence of
riggings of the partitions $\nu^{(a)}$ such that
every label $x$ of a part of $\nu^{(a)}$ of size $i$
satisfies the inequality
\begin{equation*}
  0 \le x \le p^{(a)}_i.
\end{equation*}
Alternatively, a rigging of a configuration $\nu$ may be viewed
as a double-sequence of partitions $J=(J^{(a,i)}\mid (a,i)\in\HH)$
where $J^{(a,i)}$ is a partition that has at most $m_i^{(a)}$ parts
each not exceeding $p_i^{(a)}$.
The pair $(\nu,J)$ is called a \textbf{rigged configuration}.
The set of riggings of admissible $(L,\La)$-configurations
is denoted by $\RCb(L,\La)$.
Let $(\nu,J)^{(a)}$ be the $a$-th rigged partition
of $(\nu,J)$. The \textbf{colabel} or \textbf{coquantum number} of a string $(i,x)$ 
in $(\nu,J)^{(a)}$ is defined to be $p_i^{(a)}-x$.
A string $(i,x)\in (\nu,J)^{(a)}$
is said to be \textbf{singular} if $x=p^{(a)}_i$, that is,
its label takes on the maximum value.

Using~\eqref{eq:vac}, one may easily verify that
\begin{equation}\label{eq:mp}
-p_{i-1}^{(a)}+2p_i^{(a)}-p_{i+1}^{(a)}\ge m_i^{(a-1)}-2m_i^{(a)}+m_i^{(a+1)}.
\end{equation}
This implies in particular the \textbf{convexity} condition
\begin{equation}\label{eq:convex}
p_i^{(a)}\ge \frac{1}{2}(p_{i-1}^{(a)}+p_{i+1}^{(a)}) \qquad \text{if $m_i^{(a)}=0$.}
\end{equation}

The set of rigged configurations is endowed with a natural
statistic $\cc$ called \textbf{cocharge}. For a configuration 
$\nu\in\Confb(L,\La)$ define
\begin{equation*}
\cc(\nu)=\frac{1}{2} \sum_{(a,j),(b,k)\in \HH} (\alpha_a|\alpha_b) 
 \min(j,k) m_j^{(a)} m_k^{(b)}.
\end{equation*}
For a rigged configuration $(\nu,J)\in\RCb(L,\La)$ set
\begin{equation} \label{eq:RC charge}
  \cc(\nu,J)=\cc(\nu)+\sum_{(a,i)\in\HH} |J^{(a,i)}|,
\end{equation}
where $|J^{(a,i)}|$ is the size of partition $J^{(a,i)}$.

As mentioned in the introduction, rigged configurations correspond
to highest weight crystal elements. Let $B^{r,s}$ be a Kirillov--Reshetikhin
crystal for $(r,s)\in \HH$ and $B=B^{r_k,s_k}\otimes B^{r_{k-1},s_{k-1}}\otimes \cdots 
\otimes B^{r_1,s_1}$. Associate to $B$ the multiplicity array 
$L=(L_s^{(r)} \mid (r,s)\in \HH)$ where $L_s^{(r)}$ counts the number of tensor factors 
$B^{r,s}$ in $B$. Denote by
\begin{equation*}
 \Pathb(B,\La)=\{b\in B \mid \text{$\wt(b)=\La$, $e_i(b)$ undefined for all $i\in \J$}\}
\end{equation*}
the set of all highest weight elements of weight $\La$ in $B$. There is
a natural statistics defined on $B$, called energy function or more precisely
tail coenergy function $D:B\to \Z$ (see \cite[Eq. (5.1)]{SS:2005} for a precise definition).

The following theorem was proven in~\cite{KSS:2002} for type $A_{n-1}^{(1)}$ and 
general $B=B^{r_k,s_k}\otimes \cdots \otimes B^{r_1,s_1}$, in~\cite{S:2005} for
type $D_n^{(1)}$ and $B=B^{r_k,1}\otimes \cdots \otimes B^{r_1,1}$ and 
in~\cite{SS:2005} for type $D_n^{(1)}$ and $B=B^{1,s_k}\otimes \cdots \otimes B^{1,s_1}$.
\begin{theorem} \cite{KSS:2002,S:2005,SS:2005} \label{thm:bij}
For $\La$ a dominant weight, $B$ as above and $L$ the corresponding multiplicity
array, there is a bijection $\Phib:\Pathb(B,\La)\to\RCb(L,\La)$ which preserves 
the statistics, that is, $D(b)=\cc(\Phib(b))$ for all $b\in\Pathb(B,\La)$.
\end{theorem}

Defining the generating functions
\begin{equation}\label{eq:generating bar}
\begin{split}
\Xb(B,\La)&=\sum_{b\in\Pathb(B,\La)} q^{D(b)},\\
\Mb(L,\La)&=\sum_{(\nu,J)\in \RCb(L,\La)} q^{\cc(\nu,J)},
\end{split}
\end{equation}
we get the immediate corollary of Theorem~\ref{thm:bij}.
\begin{corollary} \cite{KSS:2002,S:2005,SS:2005}
Let $\La$, $B$ and $L$ as in Theorem~\ref{thm:bij}. Then
$\Xb(B,\La)=\Mb(L,\La)$.
\end{corollary}

\subsection{Crystal structure}\label{sec:def rc crystal}
In this section we introduce the set of unrestricted rigged configurations
$\RC(L)$ by defining a crystal structure generated from highest weight vectors 
given by elements in $\RCb(L)=\bigcup_{\La\in P^+} \RCb(L,\La)$ by the Kashiwara 
operators $e_a,f_a$.

\begin{definition} \label{def:crystal} Let $L$ be a multiplicity array.
Define the set of \textbf{unrestricted rigged configurations} $\RC(L)$
as the set generated from the elements in $\RCb(L)$ by the application of  
the operators $f_a,e_a$ for $a\in \J$ defined as follows:
\begin{enumerate}
\item
Define $e_a(\nu,J)$ by removing a box from a string of length $k$ in
$(\nu,J)^{(a)}$ leaving all colabels fixed and increasing the new
label by one. Here $k$ is the length of the string with the smallest
negative rigging of smallest length. If no such string exists,
$e_a(\nu,J)$ is undefined. 
\item
Define $f_a(\nu,J)$ by adding a box to a string of length $k$ in
$(\nu,J)^{(a)}$ leaving all colabels fixed and decreasing the new
label by one. Here $k$ is the length of the string with the smallest
nonpositive rigging of largest length. If no such string exists,
add a new string of length one and label -1.
If the result is not a valid unrestricted rigged configuration
$f_a(\nu,J)$ is undefined.
\end{enumerate}
\end{definition}

Let $(\nu,J)\in \RC(L)$. If $f_a$ adds a box to a string of length $k$ in
$(\nu,J)^{(a)}$, then the vacancy numbers change according to
\begin{equation}\label{eq:change in p}
p_i^{(b)} \mapsto p_i^{(b)}-(\alpha_a | \alpha_b)\chi(i>k),
\end{equation}
where $\chi(S)=1$ if the statement $S$ is true and $\chi(S)=0$ if $S$
is false. Similarly, if $e_a$ adds a box of length $k$ to $(\nu,J)^{(a)}$,
then the vacancy numbers change as
\begin{equation*}
p_i^{(b)} \mapsto p_i^{(b)}+(\alpha_a | \alpha_b)\chi(i\ge k).
\end{equation*}
Hence if $(\nu',J')=f_a(\nu,J)$ exists, it is easy to check that 
$e_a(\nu',J')=(\nu,J)$ and vice versa.
\begin{remark}\label{rem:lowest label}
Note that it follows from~\eqref{eq:change in p} that for a string $(i,m)$ in
$(\nu,J)\in\RC(L)$ we have $m\ge -i$. Hence, $e_a$ only removes a string of
length 1 if its label is $-1$, which in Definition~\ref{def:crystal} is
interpreted as increasing the label by one.
\end{remark}

We may define a weight function $\wt:\RC(L)\to P$ as
\begin{equation}\label{eq:rc weight}
\wt(\nu,J)=\sum_{(a,i)\in \HH} i(L_i^{(a)}\La_a - m_i^{(a)}\alpha_a)
\end{equation}
for $(\nu,J)\in\RC(L)$. It is clear from the definition that 
$\wt(f_a (\nu,J))=\wt(\nu,J)-\alpha_a$.
Define
\begin{equation*}
\RC(L,\La)=\{(\nu,J)\in \RC(L) \mid \wt(\nu,J)=\La\}.
\end{equation*}

\begin{example} Let $\geh$ be of type $A_2^{(1)}$. Let $\La=\La_1-\La_2$,
$L_1^{(1)}=L_3^{(1)}=L_2^{(2)}=1$ and all other $L_i^{(a)}=0$. Then
\begin{equation*}
(\nu,J)= \quad \yngrc(2,-1,1,-1) \quad \yngrc(3,-2)
\end{equation*}
is in $\RC(L,\La)$, where the parts of the rigging $J^{(a,i)}$ are written
next to the parts of length $i$ in partition $\nu^{(a)}$. We have
\begin{equation*}
f_1(\nu,J)= \yngrc(3,-2,1,-1) \quad \yngrc(3,-1) \quad \text{and} \quad
e_1(\nu,J)= \yngrc(2,1) \quad \yngrc(3,-3).
\end{equation*}
\end{example}

Before stating our main result, we need some preliminary properties of the crystal
operators.
\begin{lemma}\label{lem:varphi}
Let $(\nu,J)\in \RC(L)$. For fixed $a\in \J$, let $p=p_i^{(a)}$ be the 
vacancy number for large $i$ and let $s\le 0$ be the smallest nonpositive label
in $(\nu,J)^{(a)}$; if no such label exists set $s=0$. Then $\varphi_a(\nu,J)=p-s$.
\end{lemma}
\begin{proof}
By definition, $f_a$ adds a box to the largest string with $s$-rigging in 
$(\nu,J)^{(a)}$. Let the length of this string be $k_0$.
By the maximality of $k_0$, we have $p_j^{(a)}>s$ for all $j>k_0$ such that $j$
exists as a part in $\nu^{(a)}$. If $p=s$, then by the convexity~\eqref{eq:convex} of 
$p_j^{(a)}$ we have $m_j^{(a)}=0$ for $j>k_0$ and $p_{k_0+1}^{(a)}=s$. Under the application 
of $f_a$, the new label of the string of length $k_0+1$ would be $s-1$, but the 
vacancy number changes by $-2$. Hence $\varphi_a(\nu,J)=0=p-s$ in this case. 
If $p>s$, by the convexity~\eqref{eq:convex} of $p_i^{(a)}$ we have that $p_j^{(a)}>s$ for 
all $j>k_0$ and $f_a(\nu,J)$ is defined. 
The new label is $s-1$ and $p_j^{(a)}\mapsto p_j^{(a)}-2\ge s-1$ for all $j>k_0$.
By induction on $m\ge 1$, $f_a$ adds a box to the largest string with rigging $s-m$
of $f_a^m(\nu,J)^{(a)}$. Let $k_m>k_{m-1}$ be the length of this string.
For large $i$, the vacancy number is $p-2m$. Suppose that $p>s+m$.
Again by the maximality of $k_m$ and the convexity of the vacancy numbers~\eqref{eq:convex}, 
the vacancy numbers of $f_a^m(\nu,J)$ satisfy $p_j^{(a)}>s-m$ for all $j>k_m$. 
The new label of the added box is $s-m-1$ and compared to the original
vacancy number $p_j^{(a)}\mapsto p_j^{(a)}-2(m+1)\ge s-m-1$ for all $j>k_m$.
If $p=s+m$, then the new label becomes $s-m-1$, but the vacancy number
is $p-2m-2=s-m-2$. Hence $f_a^{m+1}(\nu,J)$ is not defined in this case.
This proves $\varphi_a(\nu,J)=p-s$.
\end{proof}

\begin{theorem} \label{thm:rc crystal}
Let $\geh$ be of simply-laced type. For $(\nub,\Jb)\in \RCb(L,\La)$, let $X_{(\nub,\Jb)}$ 
be the graph generated by $(\nub,\Jb)$ and $e_a,f_a$ for $a\in \J$. Then $X_{(\nub,\Jb)}$ 
is isomorphic to the crystal graph $B(\La)$.
\end{theorem}

\begin{proof}
Let $A=[A_{ab}]$ be a Cartan matrix of simply-laced type
and $\La=\sum_{a\in\J} \mu_a \La_a$.
By Theorem~\ref{thm:stembridge} it suffices to check that 
the graph $X_{(\nub,\Jb)}$ generated by the maximal element $(\nub,\Jb)$ and 
operators $e_a,f_a$ as defined in Definition~\ref{def:crystal} is $A$-regular 
and that $\varphi_a(\nub,\Jb)=\mu_a$ for all $a\in \J$.

The claim that $\varphi_a(\nub,\Jb)=\mu_a$ for all $a\in \J$ follows from 
Lemma~\ref{lem:varphi}.
Combining~\eqref{eq:vac} and~\eqref{eq:conf} we find that $p_i^{(a)}=\mu_a$
for large $i$. Note that since $(\nub,\Jb)\in \RCb(L,\La)$, it does not have
any negative riggings, so that $s=0$. Hence by Lemma~\ref{lem:varphi},
$\varphi_a(\nub,\Jb)=\mu_a$.

Next we check that $X_{(\nub,\Jb)}$ is $A$-regular. 
Let $(\nu,J)\in X_{(\nub,\Jb)}$. By Lemma~\eqref{lem:varphi}, $\varphi_a(\nu,J)$ is finite.
This proves (P1). (P2) is clear from Definition~\ref{def:crystal}.

To prove (P3) and (P4) we show that one of the following conditions hold
\begin{equation*}
(A_{ab},\Delta_a\delta_b(\nu,J),\Delta_a\vp_b(\nu,J))=(0,0,0),(-1,-1,0),(-1,0,-1).
\end{equation*}
It is clear from the definitions that the operators $[f_a,f_b]=[e_a,e_b]=0$
commute when $A_{ab}=0$, so that $\Delta_a\delta_b(\nu,J)=\Delta_a\vp_b(\nu,J)=0$
in this case. Hence assume that $A_{ab}=-1$.
Let $k_a^e$ be the length of the string in $(\nu,J)^{(a)}$ selected by $e_a$.
Let $k_b^f$ be the length of the string in $(\nu,J)^{(b)}$ selected by $f_b$. 
Under $e_a$ the vacancy number changes according to
\begin{equation}\label{eq:p change}
p_i^{(b)} \mapsto p_i^{(b)}-\chi(i\ge k_a^e).
\end{equation}
Therefore by Lemma~\eqref{lem:varphi} we have
\begin{equation}\label{eq:Delta}
\Delta_a\vp_b(\nu,J)=\begin{cases} 0 & \text{if $k_a^e\le k_b^f$,}\\
-1 & \text{if $k_a^e>k_b^f$.} \end{cases}
\end{equation}
Similarly, it follows from~\eqref{eq:p change} that
\begin{equation}\label{eq:Delta d}
\Delta_a\delta_b(\nu,J)= \begin{cases} -1 & \text{if $k_a^e\le k_b^f$,}\\
0 & \text{if $k_a^e>k_b^f$.} \end{cases}
\end{equation}
(Note that by Remark~\ref{rem:lowest label} the labels $s$ of strings of length 
$i$ in $(\nu,J)^{(b)}$ satisfy $s\ge -i$ for all $(\nu,J)\in X_{(\nub,\Jb)}$. 
Hence $e_a(\nu,J)$ does not exist if this condition does not hold for $e_a(\nu,J)$).
This proves (P3) and (P4).

For (P5) assume that $\Delta_a \delta_b(\nu,J)=0$. As before, if $A_{ab}=0$, we
have $[e_a,e_b]=[f_a,f_b]=0$ and $\nabla_b\varphi_a(\nu',J')=0$ for any $(\nu',J')$,
hence in particular for $(\nu',J')=e_a e_b(\nu,J)$. Therefore assume that
$A_{ab}=-1$. By~\eqref{eq:Delta d}, $\Delta_a \delta_b(\nu,J)=0$ implies $k_a^e>k_b^f$.
An explicit calculation yields that $(\nu',J')=e_a e_b(\nu,J)=e_b e_a(\nu,J)$.
Note that 
\begin{equation*}
\begin{split}
\nabla_b\varphi_a(\nu',J')&=\varphi_a(e_a e_b (\nu,J))-
\varphi_a(f_b e_b e_a(\nu,J))\\
&=\varphi_a(e_b(\nu,J))+1-\varphi_a(\nu,J)-1=\Delta_b\varphi_a(\nu,J).
\end{split}
\end{equation*}
Since $k_b^e\le k_b^f<k_a^e\le k_a^f$ it follows from~\eqref{eq:Delta} that 
$\nabla_b\varphi_a(\nu',J')=\Delta_b\varphi_a(\nu,J)=0$.

For (P6) assume that $\Delta_a\delta_b(\nu,J)=\Delta_b\delta_a(\nu,J)=-1$.
In this case $k_a^e\le k_b^f$ and $k_b^e\le k_a^f$, and by definition
$k_a^e\le k_a^f$, $k_b^e\le k_b^f$. It can be checked explicitly that
$(\nu',J'):=e_a e_b^2 e_a(\nu,J)=e_b e_a^2 e_b(\nu,J)$ in this case.
Also, 
\begin{multline*}
\nabla_a\varphi_b(\nu',J')=\varphi_b(e_b e_a^2 e_b(\nu,J))-
\varphi_b(f_a e_a e_b^2 e_a(\nu,J))\\
=\varphi_b(e_a^2 e_b(\nu,J))+1-\varphi_b(e_a(\nu,J))-2
=\varphi_b(e_a^2 e_b(\nu,J))-\varphi_b(e_a(\nu,J))-1.
\end{multline*}
It can be shown explicitly that $\varphi_b(e_a^2 e_b(\nu,J))=\varphi_b(e_a(\nu,J))
=\varphi_b(\nu,J)$, which implies that $\nabla_a\varphi_b(\nu',J')=-1$ and similarly
with $a$ and $b$ interchanged.

(P5') and (P6') can be proved analogously.
\end{proof}

\begin{example} Consider the crystal $B(\La_1+\La_2)$ of type $A_2$
in $B=(B^{1,1})^{\otimes 3}$. Here is the crystal graph in the usual
labeling and the rigged configuration labeling:

\begin{picture}(100,200)(0,0)
\put(45,180){121}
\put(10,140){221}
\put(10,100){231}
\put(10,60){331}
\put(45,20){332}
\put(80,140){131}
\put(80,100){132}
\put(80,60){232}
\LongArrow(50,178)(25,152)
\LongArrow(60,178)(85,152)
\LongArrow(20,138)(20,112)
\LongArrow(20,98)(20,72)
\LongArrow(90,138)(90,112)
\LongArrow(90,98)(90,72)
\LongArrow(25,58)(50,32)
\LongArrow(85,58)(60,32)
\PText(30,165)(0)[b]{1}
\PText(15,120)(0)[b]{2}
\PText(15,80)(0)[b]{2}
\PText(30,40)(0)[b]{1}
\PText(95,120)(0)[b]{1}
\PText(95,80)(0)[b]{1}
\PText(80,40)(0)[b]{2}
\PText(80,165)(0)[b]{2}
\end{picture}
\hspace{1.5cm}
\begin{picture}(100,200)(0,0)
\put(40,185){{\tiny $\yngrc(1,0) \; \emptyset$}}
\put(0,145){{\tiny $\yngrc(2,-1) \; \emptyset$}}
\put(0,105){{\tiny $\yngrc(2,0) \; \yngrc(1,-1)$}}
\put(0,65){{\tiny $\yngrc(2,1) \; \yngrc(2,-2)$}}
\put(30,17){{\tiny $\yngrc(2,-1,1,-1) \; \yngrc(2,-1)$}}
\put(80,145){{\tiny $\yngrc(1,1) \; \yngrc(1,-1)$}}
\put(80,103){{\tiny $\yngrc(1,-1,1,-1) \; \yngrc(1,0)$}}
\put(80,63){{\tiny $\yngrc(2,-2,1,-1) \; \yngrc(1,0)$}}
\LongArrow(50,176)(25,152)
\LongArrow(60,176)(85,152)
\LongArrow(20,138)(20,112)
\LongArrow(20,97)(20,72)
\LongArrow(90,138)(90,117)
\LongArrow(90,94)(90,77)
\LongArrow(25,56)(50,32)
\LongArrow(83,53)(60,32)
\PText(30,165)(0)[b]{1}
\PText(15,120)(0)[b]{2}
\PText(15,80)(0)[b]{2}
\PText(30,40)(0)[b]{1}
\PText(95,120)(0)[b]{1}
\PText(95,80)(0)[b]{1}
\PText(80,40)(0)[b]{2}
\PText(80,165)(0)[b]{2}
\end{picture}
\end{example}

\begin{theorem}\label{thm:charge}
Let $X_{(\nub,\Jb)}$ be as in Theorem~\ref{thm:rc crystal}. The cocharge
$\cc$ as defined in~\eqref{eq:RC charge} is constant on $X_{(\nub,\Jb)}$.
\end{theorem}
\begin{proof}
Let $(\nu,J)\in X_{(\nub,\Jb)}$ such that $f_a(\nu,J)$ is defined. It is easy to check 
that adding a box to a string of length $k$ in $(\nu,J)^{(a)}$ changes the cocharge by
\begin{equation*}
\cc(\nu) \mapsto \cc(\nu)+1+\sum_{b,i}(\alpha_a | \alpha_b)\chi(i>k)m_i^{(b)}.
\end{equation*}
Since $f_a$ changes the label of the new string by $-1$ and leaves the 
colabels of all other strings unchanged, it is clear comparing with~\eqref{eq:change in p}
that $f_a$ does not change the total cocharge, that is $\cc(\nu,J)=\cc(f_a(\nu,J))$.
\end{proof}

For $B=B^{r_k,s_k}\otimes \cdots \otimes B^{r_1,s_1}$ and $\La\in P$ let
\begin{equation*}
\Path(B,\La)=\{ b\in B \mid \wt(b)=\La\}.
\end{equation*}

\begin{theorem} \label{thm:bij new}
Let $\La\in P$, $B$ be as in Theorem~\ref{thm:bij} and $L$ the corresponding
multiplicity array. Then there is a bijection $\Phi:\Path(B,\La)\to\RC(L,\La)$ 
which preserves the statistics, that is, $D(b)=\cc(\Phi(b))$ for all $b\in\Path(B,\La)$.
\end{theorem}
\begin{proof}
By Theorem~\ref{thm:bij} there is such a bijection for the maximal elements
$b\in \Pathb(B)$. By Theorems~\ref{thm:rc crystal} and~\ref{thm:charge} this 
extends to all of $\Path(B,\La)$.
\end{proof}

Extending the definitions of~\eqref{eq:generating bar} to
\begin{equation}\label{eq:XM def}
\begin{split}
X(B,\La)&=\sum_{b\in\Path(B,\La)} q^{D(b)},\\
M(L,\La)&=\sum_{(\nu,J)\in \RC(L,\La)} q^{\cc(\nu,J)},
\end{split}
\end{equation}
we obtain the corollary:
\begin{corollary}\label{cor:X=M}
With all hypotheses of Theorem~\ref{thm:bij new}, we have $X(B,\La)=M(L,\La)$.
\end{corollary}

\section{Unrestricted rigged configurations for type $A_{n-1}^{(1)}$}
\label{sec:type A}

In this section we give an explicit description of the elements in
$\RC(L,\La)$ for type $A_{n-1}^{(1)}$. Generally speaking, the elements
are rigged configurations where the labels lie between the vacancy number
and certain lower bounds defined explicitly (Definition~\ref{def:extended} and
Theorem~\ref{thm:ext=unres}). We use this in
Section~\ref{sec:fermi} to write down an explicit fermionic formula
for the unrestricted configuration sum $X(B,\La)$. Section~\ref{sec:affine}
is devoted to the affine crystal structure of $\RC(L,\La)$.

\subsection{Characterization of unrestricted rigged configurations}\label{sec:lower bound}
Let $L=(L_i^{(a)} \mid (a,i)\in \HH)$ be a multiplicity array and 
$\La\in P$. Recall that the set of $(L,\La)$-configurations $\Conf(L,\La)$
is the set of all sequences of partitions $\nu=(\nu^{(a)} \mid a\in \J)$ such
that~\eqref{eq:conf} holds.
As discussed in Section~\ref{sec:def rc}, in the usual setting a rigged 
configuration $(\nu,J)\in \RCb(L,\La)$ consists of a configuration
$\nu\in \Conf(L,\La)$ together with a double sequence of partitions
$J=\{J^{(a,i)}\mid (a,i)\in\HH \}$ such that the partition
$J^{(a,i)}$ is contained in a $m_i^{(a)}\times p_i^{(a)}$ rectangle.
In particular this requires that $p_i^{(a)}\ge 0$. The unrestricted rigged
configurations $(\nu,J)\in \RC(L,\La)$ can contain labels that are negative,
that is, the lower bound on the parts in $J^{(a,i)}$ can be less than zero.

To define the lower bounds we need the following notation. 
Let $\la=(\la_1,\ldots,\la_n)$ be the $n$-tuple of nonnegative integers
corresponding to $\La$, that is $\La=\sum_{i \in \J} (\la_i-\la_{i+1})\La_i$.
In this section we use $\La$ and $\la$ interchangeably.
Let $\la'=(c_1,c_2,\ldots,c_{n-1})^t$,
where $c_k=\la_{k+1}+\la_{k+2}+\cdots+\la_n$ is the length of the $k$-th
column of $\la'$, and let $\A(\la')$ be the set of tableaux of shape $\la'$ such that 
the entries are strictly decreasing along columns, and the letters in column $k$ are from 
the set $\{1,2,\ldots,c_{k-1}\}$ with $c_0=c_1$.

\begin{example} For $n=4$ and $\la=(0,1,1,1)$, the set $\A(\la')$
consists of the following tableaux
\begin{equation*}
\young(332,22,1) \quad \young(332,21,1) \quad \young(322,21,1) \quad 
\young(331,22,1) \quad \young(331,21,1) \quad \young(321,21,1).
\end{equation*}
\end{example}

\begin{remark}\label{rem:row}
Denote by $t_{j,k}$ the entry of $t\in\A(\la')$ in row $j$ and column $k$.
Note that $c_k-j+1\le t_{j,k}\le c_{k-1}-j+1$ since the entries in column $k$
are strictly decreasing and lie in the set $\{1,2,\ldots, c_{k-1}\}$.
This implies $t_{j,k}\le c_{k-1}-j+1\le t_{j,k-1}$, so that the rows of $t$
are weakly decreasing.
\end{remark}

Given $t\in\A(\la')$, we define the \textbf{lower bound} as
\begin{equation*}
M_i^{(a)}(t)=-\sum_{j=1}^{c_a} \chi(i\ge t_{j,a})
+\sum_{j=1}^{c_{a+1}} \chi(i\ge t_{j,a+1}),
\end{equation*}
where recall that $\chi(S)=1$ if the the statement $S$ is true and $\chi(S)=0$ otherwise.

Let $M,p,m\in \Z$ such that $m\ge 0$.
A $(M,p,m)$-quasipartition $\mu$ is a tuple of integers $\mu=(\mu_1,\mu_2,\ldots,\mu_m)$
such that $M\le \mu_m\le \mu_{m-1}\le \cdots\le \mu_1\le p$. Each $\mu_i$ is called
a part of $\mu$. Note that for $M=0$ this would be a partition with at most $m$ parts each 
not exceeding $p$.

\begin{definition}\label{def:extended}
An \textbf{extended rigged configuration} $(\nu,J)$ is a configuration
$\nu\in\Conf(L,\la)$ together with a sequence $J=\{J^{(a,i)}\mid (a,i)\in\HH\}$
where $J^{(a,i)}$ is a $(M_i^{(a)}(t),p_i^{(a)},m_i^{(a)})$-quasipartition
for some $t\in \A(\la')$. Denote the set of all extended rigged configurations
corresponding to $(L,\la)$ by $\RCt(L,\la)$.
\end{definition}

\begin{example}
Let $n=4$, $\la=(2,2,1,1)$, $L_1^{(1)}=6$ and all other $L_i^{(a)}=0$. Then
\begin{equation*}
(\nu,J) \;=\; \yngrc(3,-2,1,0) \quad \yngrc(2,0) \quad \yngrc(1,-1)
\end{equation*}
is an extended rigged configuration in $\RCt(L,\la)$, where we have written
the parts of $J^{(a,i)}$ next to the parts of length $i$ in partition $\nu^{(a)}$.
To see that the riggings form quasipartitions, let us write
the vacancy numbers $p_i^{(a)}$ next to the parts of length $i$ in partition $\nu^{(a)}$:
\begin{equation*}
\yngrc(3,0,1,3) \quad \yngrc(2,0) \quad \yngrc(1,-1).
\end{equation*}
This shows that the labels are indeed all weakly below the vacancy numbers. For
\begin{equation*}
\young(441,33,2,1) \in \A(\la')
\end{equation*}
we get the lower bounds
\begin{equation*}
\yngrc(3,-2,1,-1) \quad \yngrc(2,0) \quad \yngrc(1,-1),
\end{equation*}
which are less or equal to the riggings in $(\nu,J)$.
\end{example}

\begin{remark}\mbox{}
\begin{enumerate}
\item
Note that Definition~\ref{def:extended} is similar to the definition of level-restricted 
rigged configurations~\cite[Definition 5.5]{SS:2001}. Whereas for level-restricted
rigged configurations the vacancy number had to be modified according to tableaux
in a certain set, here the lower bounds are modified.
\item 
For type $A_1$ we have $\la=(\la_1,\la_2)$ so that $\A=\{ t\}$
contains just the single tableau
\begin{equation*}
t=\begin{array}{|c|} \hline \la_2\\ \hline \la_2-1\\ \hline \vdots\\ \hline 1\\
\hline \end{array}.
\end{equation*}
In this case $M_i(t)=-\sum_{j=1}^{\la_2} \chi(i\ge t_{j,1})=-i$. This agrees
with the findings of~\cite{T:2005}.
\end{enumerate}
\end{remark}

The next theorem shows that the set of unrestricted rigged configurations $\RC(L,\la)$
of type $A_{n-1}^{(1)}$ defined in terms of the crystal structure in 
Section~\ref{sec:def rc crystal} is equal to the set of extended rigged configurations
$\RCt(L,\la)$ of Definition~\ref{def:extended}.

\begin{theorem} \label{thm:ext=unres}
We have $\RC(L,\la)=\RCt(L,\la)$.
\end{theorem}
\begin{proof}
Denote by $X_{(\nub,\Jb)}$ the graph with maximal element $(\nub,\Jb)\in \RCb(L)$ 
generated by $f_a,e_a$ for $a\in \J$. By definition
\begin{equation*}
\bigcup_{(\nub,\Jb)\in\RCb(L)} X_{(\nub,\Jb)} = \RC(L).
\end{equation*}
We claim that $\RC(L)=\RCt(L)$. The statement $\RC(L,\la)=\RCt(L,\la)$ then follows
since the weight function is defined in the same way on both sets.

Let $(\nub,\Jb)\in\RCb(L,\la)$. Then $(\nub,\Jb)$ is admissible with respect to
$t\in\A(\la')$ where column $k$ of $t$ is filled with the letters $1,2,\ldots,c_k$.
Hence $(\nub,\Jb)\in \RCt(L)$.
Now suppose by induction that $(\nu,J)\in \bigcup_{(\nub,\Jb)\in\RCb(L)} X_{(\nub,\Jb)}$
is admissible with respect to $t\in\A(\la')$. We claim that $(\nu',J')=f_a(\nu,J)$ 
is admissible with respect to some $t'\in\A(\la'')$ where $\la''$ is obtained from 
$\la'$ by adding a box to column $a$. (Note that $f_a(\nu,J)=0$ if $\la_a=0$ or 
equivalently $c_{a-1}=c_a$).
Let $k$ be the length of the string in $(\nu,J)^{(a)}$ selected by $f_a$ (see
Definition~\ref{def:crystal}).
Let $r>k$ be minimal such that $r\not\in t_{\cdot,a}$, where $t_{\cdot,a}$ denotes
column $a$ of $t$. Similarly, let $s>k$ be minimal such that $s\in t_{\cdot,a+1}$. 
Then $t'$ is obtained from $t$ by adding $r$ to column $a$, and by removing $s$ from
column $a+1$ and adding $c_a+1$ to column $a+1$ in such a way that the columns are 
still strictly decreasing. Note that $t'$ is by construction strictly decreasing in
columns and has the property that the elements in column $b$ lie in the set
$\{1,2,\ldots,c'_{b-1}\}$ where $c'_b=c_b+\delta_{a,b}$ is the length of column 
$b$ in $t'$. Hence $t'\in\A(\la'')$.

To see that $(\nu',J')=f_a(\nu,J)$ is admissible with respect to $t'$, note that 
strings in $(\nu,J)^{(a-1)}$ and $(\nu,J)^{(a+1)}$ change by 
$(j,x)\mapsto (j,x+\chi(j>k))$, and strings in $(\nu,J)^{(a)}$ change by
$(j,x)\mapsto (j,x-2\chi(j>k))$. In addition to this there is a new string $(k+1,-m-1)$ 
in $(\nu',J')^{(a)}$ where $-m$ is the smallest label in $(\nu,J)^{(a)}$.
Since column $a$ of $t'$ contains an additional entry greater than $k$ and
in column $a+1$ an entry greater than $k$ was increased, $t'$ certainly provides
valid lower bounds for $(\nu',J')^{(a\pm 1)}$. Note that
\begin{equation}\label{eq:M change}
M_{j+1}^{(a)}(t)=M_j^{(a)}(t)-\chi(j+1\in t_{\cdot,a})+\chi(j+1\in t_{\cdot,a+1}).
\end{equation}
Since by definition of $f_a$, $k$ is largest such that there is a string of this
length with label $-m$, it is not hard to check that $t'$ gives proper lower
bounds for $(\nu',J')$. This shows that
\begin{equation*}
\bigcup_{(\nub,\Jb)\in\RCb(L)} X_{(\nub,\Jb)} \subset \RCt(L).
\end{equation*}

To prove the reverse inclusion, suppose $(\nu,J)\in \RCt(L)$. Let $t\in \A(\la')$ be such
that $(\nu,J)$ is admissible with respect to $t$. If $(\nu,J)\in \RCb(L)$, then certainly
$(\nu,J)\in \bigcup_{(\nub,\Jb)\in\RCb(L)} X_{(\nub,\Jb)}$. If $(\nu,J)\not\in \RCb(L)$,
there must be at least one negative rigging. Suppose this occurs in $(\nu,J)^{(a)}$.
Then $e_a(\nu,J)$ exists. To see this note that all colabels remain fixed, so that
all labels are still weakly below the vacancy number. The string $(k,-m)$ in 
$(\nu,J)^{(a)}$ selected by $e_a$ becomes $(k-1,-m+1)$. Since by the definition
of $k$, $k$ is smallest such that its label is $-m<0$, all labels of strings of length
less than $k$ are strictly bigger than $-m$. Hence $p_j^{(a)}>-m$ for all $j<k$
such that $j$ appears as a part. By the convexity property~\eqref{eq:convex} of $p_j^{(a)}$, 
this is true for all $j<k$. Hence $e_a(\nu,J)$ exists. 

Next we need to show that $(\nu',J')=e_a(\nu,J)\in \RCt(L)$. Let $r\le k$ be maximal 
such that $r\in t_{\cdot,a}$ and let $s\le k$ be maximal such that 
$s\not\in t_{\cdot,a+1}$. Note that $r$ and $s$ must exist, since the
rigging of the string of length $k$ in $(\nu,J)^{(a)}$ is negative so that
$M_k^{(a)}(t)<0$. But this implies that $\#\{j\in t_{\cdot,a} \mid j\le k\}>
\#\{j\in t_{\cdot,a+1} \mid j\le k\}$.
Then define $t'$ by removing $r$ from column $a$ of $t$ and
changing the largest element in column $a+1$ to $s$. 
By similar arguments as for the previous case, $e_a(\nu,J)$ is admissible
with respect to $t'$.
\end{proof}

\begin{remark}
For type $D_n^{(1)}$, a simple characterization in terms of lower bounds
for the parts of a configuration $\nu\in C(L)$ does not seem to exist.
For example take $B=B^{2,1}$ of type $D_4^{(1)}$ so that $L_1^{(2)}=1$ and all 
other $L_i^{(a)}=0$. Then the unrestricted rigged configurations
\begin{equation*}
\yngrc(1,0) \quad \yngrc(1,0,1,0) \quad \yngrc(1,0) \quad \yngrc(1,0)
\quad \text{and}\quad
\yngrc(1,0) \quad \yngrc(1,0,1,-1) \quad \yngrc(1,0) \quad \yngrc(1,0),
\end{equation*}
which correspond to the crystal elements $\begin{array}{|c|} \hline 1\\ \hline
\overline{1} \\ \hline \end{array}$ and $\begin{array}{|c|} \hline 3\\ \hline
\overline{3}\\ \hline \end{array}$ respectively, occur in $\RC(L)$, but
\begin{equation*}
\yngrc(1,0) \quad \yngrc(1,-1,1,-1) \quad \yngrc(1,0) \quad \yngrc(1,0)
\end{equation*}
on the other hand does not appear.
\end{remark}

\subsection{Fermionic formula} \label{sec:fermi}
With the explicit characterization of the unrestricted rigged configurations
of Section~\ref{sec:lower bound}, it is possible to derive an explicit formula 
for the polynomials $M(L,\la)$ of~\eqref{eq:XM def}.

Let $\SA(\la')$ be the set of all nonempty subsets of $\A(\la')$ and set
\begin{equation*}
M_i^{(a)}(S)=\max\{M_i^{(a)}(t) \mid t\in S\} \qquad \text{for $S\in\SA(\la')$.}
\end{equation*}
By inclusion-exclusion the set of all allowed riggings for a given $\nu\in\Conf(L,\la)$ is
\begin{equation*}
\bigcup_{S\in\SA(\la')} (-1)^{|S|+1} \{J\mid \text{$J^{(a,i)}$ is a 
 $(M_i^{(a)}(S),p_i^{(a)},m_i^{(a)})$-quasipartition}\}.
\end{equation*}
The $q$-binomial coefficient $\qbin{m+p}{m}$, defined as
\begin{equation*}
\qbin{m+p}{m}=\frac{(q)_{m+p}}{(q)_m(q)_p},
\end{equation*}
where $(q)_n=(1-q)(1-q^2)\cdots(1-q^n)$, is the generating function of partitions
with at most $m$ parts each not exceeding $p$. Hence the polynomial $M(L,\la)$
may be rewritten as 
\begin{multline} \label{eq:fermi}
M(L,\la)=\sum_{S\in\SA(\la')} (-1)^{|S|+1} \sum_{\nu\in\Conf(L,\la)}
q^{\cc(\nu)+\sum_{(a,i)\in\HH} m_i^{(a)}M_i^{(a)}(S)}\\
\times \prod_{(a,i)\in\HH} \qbin{m_i^{(a)}+p_i^{(a)}-M_i^{(a)}(S)}{m_i^{(a)}}
\end{multline}
called \textbf{fermionic formula}. By Corollary~\ref{cor:X=M} this is also a formula
for the unrestricted configuration sum $X(B,\la)$. As mentioned in the introduction,
this formula is different from the fermionic formulas of~\cite{HKKOTY:1999,Kir:2000}
which exist in the special case when $L$ is the multiplicity array of
$B=B^{1,s_k}\otimes \cdots \otimes B^{1,s_1}$ or $B=B^{r_k,1}\otimes \cdots 
\otimes B^{r_1,1}$.

\subsection{The Kashiwara operators $e_0$ and $f_0$}\label{sec:affine}
The Kirillov--Reshetikhin crystals $B^{r,s}$ are affine crystals and admit
the Kashiwara operators $e_0$ and $f_0$. It was shown in~\cite{Sh:2002} that
for type $A_{n-1}^{(1)}$ they can be defined in terms of the
\textbf{promotion operator} $\pr$ as
\begin{equation*}
e_0=\pr^{-1}\circ e_1\circ \pr \quad \text{and} \quad 
f_0=\pr^{-1}\circ f_1\circ \pr.
\end{equation*}
The promotion operator is a bijection $\pr:B\to B$ such that the following
diagram commutes for all $a\in I$ 
\begin{equation}\label{eq:pr com}
\begin{CD}
B@>{\pr}>> B \\
@V{f_a}VV @VV{f_{a+1}}V \\
B @>>{\pr}> B
\end{CD}
\end{equation}
and such that for every $b\in B$ the weight is rotated
\begin{equation}\label{eq:weight rotation}
\inner{h_{a+1}}{\wt(pr(b))}=\inner{h_a}{\wt(b)}.
\end{equation}
Here subscripts are taken modulo $n$. 

We are now going to define the promotion operator on unrestricted rigged configurations.
\begin{definition}\label{def:rc pr}
Let $(\nu,J)\in\RC(L,\la)$. Then $\pr(\nu,J)$ is obtained as follows:
\begin{enumerate}
\item Set $(\nu',J')=f_1^{\la_1} f_2^{\la_2} \cdots f_n^{\la_n}(\nu,J)$ where
$f_n$ acts on $(\nu,J)^{(n)}=\emptyset$.
\item Apply the following algorithm $\rho$ to $(\nu',J')$ $\la_n$ times: Find the smallest
singular string in $(\nu',J')^{(n)}$. Let the length be $\ell^{(n)}$.
Repeatedly find the smallest singular string in $(\nu',J')^{(k)}$ of length
$\ell^{(k)}\ge \ell^{(k+1)}$ for all $1\le k<n$. Shorten the selected strings
by one and make them singular again.
\end{enumerate}
\end{definition}

\begin{example}
Let $B=B^{2,2}$, $L$ the corresponding multiplicity array and $\la=(1,0,1,2)$. Then
\begin{equation*}
(\nu,J)=\quad \yngrc(1,0) \quad \yngrc(2,-1,1,-1) \quad \yngrc(2,-1) \quad
\in \RC(L,\la)
\end{equation*}
corresponds to the tableau $b=\young(13,44)\in \Path(B,\la)$. After step (1)
of Definition~\ref{def:rc pr} we have
\begin{align*}
(\nu',J')&=\quad \yngrc(2,-1) \quad \yngrc(2,1,1,0) \quad \yngrc(2,-1,1,-1) \quad
\yngrc(2,-1).\\
\intertext{Then applying step (2) yields}
\pr(\nu,J)&=\quad \emptyset \quad \yngrc(1,0) \quad \yngrc(1,-1)
\end{align*}
which corresponds to the tableau $\pr(b)=\young(11,24)$.
\end{example}

\begin{lemma}\label{lem:pr}
The map $\pr$ of Definition~\ref{def:rc pr} is well-defined and 
satisfies~\eqref{eq:pr com} for $1\le a\le n-2$ and~\eqref{eq:weight rotation} 
for $0\le a\le n-1$.
\end{lemma}
\begin{proof}
To prove that $\pr$ is well-defined we need to show that singular strings of length 
$\ell^{(k)}$ exist in $(\nu',J')^{(k)}=f_1^{\la_1} f_2^{\la_2} \cdots f_n^{\la_n}(\nu,J)$,
where $\ell^{(1)}\ge \ell^{(2)}\ge \cdots \ge \ell^{(n)}$.
For a given $1\le a\le n$, set $(\tilde{\nu},\tilde{J})=f_a^{\la_a} f_{a+1}^{\la_{a+1}}
\cdots f_n^{\la_n}(\nu,J)$. By definition $\varphi_a(\tilde{\nu},\tilde{J})=0$.
Hence by Lemma~\ref{lem:varphi} $p=s$ where $p=p_i^{(a)}$ for large $i$ and $s$ is the
smallest nonpositive label in $(\tilde{\nu},\tilde{J})^{(a)}$. Let $\ell$ be the
length of the largest part in $(\tilde{\nu},\tilde{J})^{(a)}$. Suppose that
$(\tilde{\nu},\tilde{J})^{(a-1)}$ or $(\tilde{\nu},\tilde{J})^{(a+1)}$ has a part
of length bigger than $\ell$. Then by the definition of the vacancy number,
$p_\ell^{(a)}<p$. But this contradicts the fact that $s=p$ is the smallest label
in $(\tilde{\nu},\tilde{J})^{(a)}$. Hence the parts of $(\tilde{\nu},\tilde{J})^{(a\pm 1)}$
cannot exceed $\ell$ and the string $(\ell,p)$ in $(\tilde{\nu},\tilde{J})^{(a)}$
is singular. Since the application of $f_1^{\la_1} \cdots f_{a-1}^{\la_{a-1}}$ does
not change the colabels in the $a$-th rigged partition, the largest string remains
singular. Note that the above argument also shows that the longest parts in ${\nu'}^{(a)}$
decrease with $a$. Hence there exist singular strings in $(\nu',J')$ such
that $\ell^{(k)}\ge \ell^{(k+1)}$ and $\pr$ is well-defined.

Next we show that $\pr$ satisfies~\eqref{eq:weight rotation}. Let $(\nu,J)\in\RC(L,\la)$
so that $\wt(\nu,J)=\la=(\la_1,\ldots,\la_n)$. After step (1) of 
Definition~\ref{def:rc pr}, we have $|\nu^{(a)}|\mapsto |\nu^{(a)}|+\la_a$.
Hence by~\eqref{eq:rc weight}, noting that $\sum_i i m_i^{(a)}=|\nu^{(a)}|$, we obtain
\begin{equation*}
\begin{split}
\wt(\nu',J')&=\wt(\nu,J)-\sum_{a\in\J} \la_a \alpha_a\\
&=(\la_1,\ldots,\la_n,0)+(-\la_1,\la_1-\la_2,\ldots,\la_{n-1}-\la_n,\la_n)
=(0,\la_1,\ldots,\la_n)
\end{split}
\end{equation*}
where in the second line we added an $(n+1)$-th component to the weight.
In step (2) of Definition~\ref{def:rc pr} the size of the $a$-th partition
changes as $|\nu^{(a)}|\mapsto |\nu^{(a)}|-\la_n$. Hence
\begin{equation*}
\begin{split}
\wt(\pr(\nu,J))&=\wt(\nu',J')-\sum_{a\in\J} \la_n \alpha_a\\
&=(0,\la_1,\ldots,\la_n)-(-\la_n,0,\ldots,0,\la_n)=(\la_n,\la_1,\la_2,\ldots,\la_{n-1},0).
\end{split}
\end{equation*}
Dropping the last component (which we only added for the intermediate calculation)
we obtain~\eqref{eq:weight rotation}.

It remains to prove~\eqref{eq:pr com}. We treat step (1) and step (2) in the
Definition~\ref{def:rc pr} of $\pr$ separately. 
Note that $f_a$ and $f_b$ commute as long as $b\neq a\pm 1$.
Hence for step (1) it suffices to show that for $(\nu,J)\in\RC(L,\la)$ with $\la_{a+2}=0$
we have $f_{a+1} f_a^{\la_a} f_{a+1}^{\la_{a+1}} (\nu,J)=f_a^{\la_a-1} f_{a+1}^{\la_{a+1}+1} 
f_a(\nu,J)$. Note that it is not hard to check that Lemma~\ref{lem:varphi} implies
\begin{equation*}
\nabla_a \varphi_{a+1}(\nu,J)=\begin{cases}
 0 & \text{if $k_a^f\le k_{a+1}^e$,}\\
 -1 & \text{if $k_a^f>k_{a+1}^e$,}
\end{cases}
\end{equation*}
where $k_a^f$ (resp. $k_{a+1}^e$) is the length of the string in $(\nu,J)^{(a)}$
(resp. $(\nu,J)^{(a+1)}$) selected by $f_a$ (resp. $e_{a+1}$).
Since in our case $\nabla_a \varphi_{a+1}(\nu,J)=-1$ we must have $k_a^f>k_{a+1}^e$.
Let $\ell_a^f$ be the length of the string selected by $f_a$ in 
$f_{a+1}^{\la_{a+1}}(\nu,J)$. If $\ell_a^f>k_{a+1}^f$, then necessarily $k_a^f=\ell_a^f$
since the application of $f_{a+1}$ only increases labels in the $a$-th rigged partition.
Note that in this case $f_a f_{a+1}(\nu,J)=f_{a+1} f_a(\nu,J)$. Hence it suffices
to repeat the analysis for $f_a f_{a+1}(\nu,J)$.
If $\ell_a^f\le k_{a+1}^f$, it can be checked explicitly that 
$f_{a+1} f_a^{\la_a} f_{a+1}^{\la_{a+1}} (\nu,J)=f_a^{\la_a-1} f_{a+1}^{\la_{a+1}+1} 
f_a(\nu,J)$.

The algorithm $\rho$ in step (2) of Definition~\ref{def:rc pr} commutes with
$f_a$ for all $2\le a\le n-1$ (assuming that both $\rho$ and $f_a$ are defined on
the rigged configuration)
\begin{equation}\label{eq:rho com}
\begin{CD}
\RC(L)@>{\rho}>> \RC(L) \\
@V{f_a}VV @VV{f_a}V \\
\RC(L) @>>{\rho}> \RC(L).
\end{CD}
\end{equation}
To see this, let $(\nu,J)\in\RC(L)$, let $(k,s)$ be the string of length $k$ and
label $s$ selected by $f_a$ in $(\nu,J)$ and denote by $\ell^{(b)}$ the length of the 
strings selected by $\rho$ in $(\nu,J)$. If $k>\ell^{(a)}$ or $k<\ell^{(a+1)}-1$, then 
$\rho$ and $f_a$ obviously commute since $f_a$ leaves all colabels of unselected strings 
fixed and $\rho$ leaves all labels of unselected strings fixed. Hence assume that 
$\ell^{(a+1)}-1\le k\le \ell^{(a)}$. Since $s$ is the smallest label occurring in
$(\nu,J)^{(a)}$, we must have $p_i^{(a)}\ge s$ for all $i$ such that $m_i^{(a)}>0$.
Because of the convexity property~\eqref{eq:convex} of the vacancy numbers, it follows 
that $p_i^{(a)}\ge s$ for all $i$. Since $(k,s)$ is the largest string in $(\nu,J)^{(a)}$
with label $s$, all strings $(i,x)$ with $i>k$ in $(\nu,J)^{(a)}$ satisfy $x>s$. 
This together with the convexity condition~\eqref{eq:convex} implies that $p_i^{(a)}\ge s+1$ 
for all $i>k$. Also, since there are no singular strings of length $k<i<\ell^{(a)}$ in 
$(\nu,J)^{(a)}$, we must have $p_i^{(a)}\ge s+2$ for $k<i<\ell^{(a)}$ and $m_i^{(a)}>0$.

Let us first assume that $k=\ell^{(a)}$. Then either $m_k^{(a)}\ge 2$, or 
$m_k^{(a)}=1$ and $p_k^{(a)}=s$. 
First consider the case that $m_k^{(a)}=1$ and $p_k^{(a)}=s$.
For $i=k$ the inequality~\eqref{eq:mp} reads
\begin{equation}\label{eq:mp in}
2\ge m_k^{(a-1)}+m_k^{(a+1)}+(p_{k-1}^{(a)}-p_k^{(a)})+(p_{k+1}^{(a)}-p_k^{(a)}).
\end{equation}
Certainly $p_{k+1}^{(a)}>p_k^{(a)}=s$ as discussed above.
Similarly, $p_{k-1}^{(a)}\ge p_k^{(a)}$ since $s$ is the smallest label in $(\nu,J)^{(a)}$.
This implies that $p_{k-1}^{(a)}=p_k^{(a)}$ or $p_{k-1}^{(a)}=p_k^{(a)}+1$.
If $\ell^{(a+1)}<\ell^{(a)}$ we must have $p_{k-1}^{(a)}=p_k^{(a)}+1$. This is clear
if $m_{k-1}^{(a)}>0$ since else $\rho$ would pick the singular string of length $k-1$ rather 
than $k$. If $m_{k-1}^{(a)}=0$ let $k'<k$ be maximal such that $m_{k'}^{(a)}>0$.
By~\eqref{eq:mp} it follows that $m_i^{(a+1)}=0$ for $k'<i<k$ and $p_i^{(a)}=s$ for
$k'\le i\le k$ if $p_k^{(a)}=p_{k-1}^{(a)}=s$. But then $\ell^{(a+1)}\le k'$ and $\rho$ 
would pick the singular string of length $k'$ in $(\nu,J)^{(a)}$ which is a contradiction 
to our assumption $\ell^{(a)}=k$. Hence $p_{k-1}^{(a)}=p_k^{(a)}+1$ if 
$\ell^{(a+1)}<\ell^{(a)}$. If $\ell^{(a+1)}=\ell^{(a)}=k$ we must have $m_k^{(a+1)}>0$ so 
that $p_{k-1}^{(a)}=p_k^{(a)}$.
Now after the application of $\rho$ there is a new string $(k-1,p_{k-1}^{(a)}
-\chi(\ell^{(a+1)}<\ell^{(a)}))=(k-1,s)$. Since this is the longest string with
label $s$, it will be picked by $f_a$ to give the string $(k,s-1)$.
If $f_a$ is applied first, there is a new string $(k+1,s-1)$.
By~\eqref{eq:mp in} and the previous findings, we have $p_{k+1}^{(a)}=s+1$.
Under $f_a$ the vacancy number changes to $p_{k+1}^{(a)}\mapsto p_{k+1}^{(a)}-2=s-1$.
Hence the string $(k+1,s-1)$ in $f_a(\nu,J)^{(a)}$ is singular and will be picked
by $\rho$. Note also that by~\eqref{eq:mp in} $m_k^{(a-1)}=0$ so that $\ell^{(a-1)}>
k=\ell^{(a)}$. This shows that $f_a$ and $\rho$ commute in this case.

If $k=\ell^{(a)}$ and $m_k^{(a)}\ge 2$ it needs to be shown that the string
of length $k$ picked by $f_a$ is still picked after the application of $\rho$.
The only case in which this might not happen is when the new string $(k-1,x)$
produced by $\rho$ has label $x<s$. Note that $p_{k-1}^{(a)}\ge s$ and
$p_{k-1}^{(a)} \mapsto p_{k-1}^{(a)}-\chi(\ell^{(a+1)}<k)=x$ under $\rho$.
Hence there is only a problem if $p_{k-1}^{(a)}=p_k^{(a)}=s$ and $\ell^{(a+1)}<k$.
Distinguishing the two cases $m_{k-1}^{(a)}>0$ and $m_{k-1}^{(a)}=0$, by very similar 
arguments as above this is not possible. Therefore $\rho$ and $f_a$ commute. 

Hence from now on we assume that $k<\ell^{(a)}=:\ell$.
Note that $p_{\ell}^{(a)}\ge s+1$ since all strings of length greater than $k$ have 
label greater than $s$. Also $p_k^{(a)}\ge s+1$ if $k\ge \ell^{(a+1)}$ since in this 
case the string $(k,s)$ cannot be singular. Convexity~\eqref{eq:convex} implies that 
$p_{\ell-1}^{(a)}\ge s+2$ and $p_{k+1}^{(a)}\ge s+2$ unless $m_i^{(a)}=0$ for all $k<i<\ell$ 
and $p_i^{(a)}=s+1$ for all $k\le i\le\ell$, or $k=\ell^{(a+1)}-1$.
If $p_{\ell-1}^{(a)}\ge s+2$ and $p_{k+1}^{(a)}\ge s+2$, $f_a$ creates a new string 
$(k+1,s-1)$ with new vacancy number $p_{k+1}^{(a)}\ge s$. Hence this string is not
singular and $\rho$ still picks the string $(\ell,p_\ell^{(a)})$. Applying $\rho$
first creates a new string $(\ell-1,p_{\ell-1}^{(a)}-1)$ with label $p_{\ell-1}^{(a)}-1\ge
s+1$. Hence $f_a$ picks the same string $(k,s)$ as before which shows that $f_a$ and
$\rho$ commute. 

Now assume that $m_i^{(a)}=0$ for $k<i<\ell$, $p_i^{(a)}=s+1$ for $k\le i\le \ell$ 
and $k\ge \ell^{(a+1)}$. In this case $f_a$ makes a new string $(k+1,s-1)$ which is 
singular since $p_{k+1}^{(a)}\mapsto p_{k+1}^{(a)}-2=s-1$. Then $\rho$ picks
this string and makes it into $(k,s)$ since under $\rho\circ f_a$ the vacancy
number changes to $p_k^{(a)}\mapsto p_k^{(a)}-1=s$. On the other hand, first applying
$\rho$ picks the string $(\ell,s+1)$ and makes it into $(\ell-1,s)$. Since this is
now the largest string with label $s$, $f_a$ picks it and makes it into
$(\ell,s-1)$ which is the same as under $\rho\circ f_a$ since recall that
$f_a$ does not change colabels but $p_\ell^{(a)}\mapsto p_\ell^{(a)}-2=s-1$ under
$f_a$. Since $m_i^{(a-1)}=0$ for all $k<i<\ell$ by~\eqref{eq:mp}, $\ell^{(a-1)}$ 
remains unchanged. Hence $f_a$ and $\rho$ commute.

Finally assume that $k=\ell^{(a+1)}-1$ and $p_{k+1}^{(a)}=s+1$. This implies that
$m_i^{(a)}=0$ for $k<i<\ell$. By~\eqref{eq:mp} with $i=k+1$ and $\ell>k+1$ we obtain
\begin{equation*}
0\ge m_{k+1}^{(a-1)}+m_{k+1}^{(a+1)}+(p_k^{(a)}-p_{k+1}^{(a)})+(p_{k+2}^{(a)}-p_{k+1}^{(a)}).
\end{equation*}
Since $\ell^{(a+1)}=k+1$ we have $m_{k+1}^{(a+1)}\ge 1$. Also $p_k^{(a)}-p_{k+1}^{(a)}
=p_k^{(a)}-s-1\ge -1$ and $p_{k+2}^{(a)}-p_{k+1}^{(a)}\ge 0$ since $p_k^{(a)}\ge s$
and $p_{k+2}^{(a)}\ge s+1$. Hence we must have $p_k^{(a)}=s$.
Now $f_a$ creates a new string $(k+1,s-1)$ which is singular since $p_{k+1}^{(a)}
\mapsto p_{k+1}^{(a)}-2=s-1$. Then $\rho$ picks this string and makes a string
$(k,s)$. If $\rho$ is applied first it transforms the string $(\ell,s+1)$ to
$(\ell-1,s)$. Then this becomes the longest string with label $s$, so that $f_a$
picks it and transforms it into $(\ell,s-1)$. Hence $f_a$ and $\rho$ commute.
If $\ell=k+1$ it is easy to show that $f_a$ and $\rho$ commute.

This concludes the proof of Lemma~\ref{lem:pr}.
\end{proof}

Lemma 7 of~\cite{Sh:2002} states that for a single Kirillov--Reshetikhin crystal
$B=B^{r,s}$ the promotion operator $\pr$ is uniquely determined by~\eqref{eq:pr com} 
for $1\le a\le n-2$ and~\eqref{eq:weight rotation} for $0\le a\le n-1$. 
Hence by Lemma~\ref{lem:pr} $\pr$ on $\RC(L)$ is indeed the correct promotion
operator when $L$ is the multiplicity array of $B=B^{r,s}$.

\begin{theorem} \label{thm:pr}
Let $L$ be the multiplicity array of $B=B^{r,s}$. Then $\pr:\RC(L)\to\RC(L)$
of Definition~\ref{def:rc pr} is the promotion operator on rigged configurations.
\end{theorem}

\begin{conjecture} \label{conj:pr}
Theorem~\ref{thm:pr} is true for any $B=B^{r_k,s_k}\otimes \cdots \otimes B^{r_1,s_1}$.
\end{conjecture}

Unfortunately, the characterization~\cite[Lemma 7]{Sh:2002} does not suffice to
define $\pr$ uniquely on tensor products $B=B^{r_k,s_k}\otimes \cdots \otimes B^{r_1,s_1}$.
In~\cite{DS:2005} a bijection $\Phi:\Path(B,\la)\to \RC(L,\la)$ is defined via a
direct algorithm. It is expected that Conjecture~\ref{conj:pr} can be proven by showing
that the following diagram commutes:
\begin{equation*}
\begin{CD}
\Path(B)@>{\Phi}>> \RC(L)\\
@V{\pr}VV @VV{\pr}V \\
\Path(B) @>>{\Phi}> \RC(L).
\end{CD}
\end{equation*}
Alternatively, an independent characterization of $\pr$ on tensor factors would give
a new, more conceptual way of defining the bijection $\Phi$ between paths and
(unrestricted) rigged configurations. A proof that the crystal operators
$f_a$ and $e_a$ commute with $\Phi$ for $a=1,2,\ldots,n-1$ is given in~\cite{DS:2005}.

\end{document}